\documentclass[a4paper,12pt]{amsart}
\usepackage{verbatim}
\usepackage{rotating} 
\usepackage[left=1.1in, right=1.1in, top=1.2in, bottom=1.2in]{geometry}
\usepackage[bookmarks,colorlinks,breaklinks,pagebackref,linktocpage]{hyperref} 
\hypersetup{linkcolor=blue,citecolor=red,filecolor=blue,urlcolor=blue}
\usepackage{mathrsfs,amsfonts,amsmath,amssymb,epsfig,amscd,xy,amsthm,bbm,mathtools} 

\usepackage{microtype}
\newtheorem{theorem}{Theorem}[section]

\newtheorem{conjecture}[theorem]{Conjecture}

\theoremstyle{plain}

\theoremstyle{remark}

\newtheorem{remark}[theorem]{Remark}
\newtheorem{example}[theorem]{Example}

\newcommand{\bbk}{\mathbbm{k}}

\newcommand{\G}{{\mathbb G}}
\newcommand{\Q}{{\mathbb Q}}

\newcommand{\N}{{\mathbb N}}
\newcommand{\K}{{\mathbb K}}
\renewcommand{\L}{{\mathbb L}}

\newcommand{\Qbar}{\overline{\Q}}

\newcommand{\bP}{{\mathbb P}}

\newcommand{\bG}{{\mathbb G}}

\newcommand{\lra}{\longrightarrow}

\author{Dragos Ghioca}
\address{
Dragos Ghioca\\
Department of Mathematics\\
University of British Columbia\\
1984 Mathematics Road\\
Vancouver, BC V6T 1Z2\\
Canada
}
\email{dghioca@math.ubc.ca}

\author{Fei Hu}
\address{
Fei Hu\\
Department of Mathematics\\
University of British Columbia\\
1984 Mathematics Road\\
Vancouver, BC V6T 1Z2\\
Canada
\endgraf
Pacific Institute for the Mathematical Sciences\\
2207 Main Mall\\
Vancouver, BC V6T 1Z4\\
Canada}
\email{fhu@math.ubc.ca}

\author{Thomas Scanlon}
\address{Thomas Scanlon\\
University of California, Berkeley\\
Mathematics Department\\
Evans Hall\\
Berkeley, CA 94720-3840\\
USA}
\email{scanlon@math.berkeley.edu}

\author{Umberto Zannier}
\address{Umberto Zannier\\
Scuola Normale Superiore\\
Piazza dei Cavalieri, 7\\
56126 Pisa\\
Italy}
\email{u.zannier@sns.it}

\thanks{The first author has been partially supported by a Discovery Grant from the National Science and Engineering Board of Canada.
The second author was partially supported by a UBC-PIMS Postdoctoral Fellowship.
The third author has been partially supported by grant DMS-1363372 of the United States National Science Foundation and a Simons Foundation Fellowship.
The authors thank the anonymous referee for making suggestions which have improved the presentation of our work.}


\begin{document}
\title{A variant of the Mordell--Lang conjecture}

\begin{abstract}
The Mordell--Lang conjecture (proven by Faltings, Vojta and McQuillan) states that the intersection of a subvariety $V$ of a semiabelian variety $G$ defined over an algebraically closed field $\bbk$ of characteristic $0$ with a finite rank subgroup $\Gamma\leq G(\bbk)$ is a finite union of cosets of subgroups of $\Gamma$. We explore a variant of this conjecture when $G=\G_a\times A$ for an abelian variety $A$ defined over $\bbk$.
\end{abstract}

\maketitle


\section{Introduction}
\label{sec:introduction}

\noindent
Throughout our paper, each subvariety is assumed to be closed. Unless otherwise noted, $\bbk$ will always denote an algebraically closed field of characteristic $0$.

Faltings \cite{Fal94} proved the Mordell--Lang conjecture, thus showing that any subvariety of an abelian variety $A$ defined over $\bbk$ intersects a finitely generated subgroup $\Gamma\leq A(\bbk)$ in a finite union of cosets of subgroups of $\Gamma$. Vojta \cite{Voj96} proved that Faltings' result holds when we replace the abelian variety $A$ by extensions $G$ of $A$ by an algebraic torus, i.e., when $G$ is a semiabelian variety. Then McQuillan \cite{McQuillan} extended further Vojta's theorem by proving that the conclusion holds when we replace $\Gamma$ by any finite rank subgroup of a semiabelian variety. It is natural to ask whether a variant of the Mordell--Lang conjecture holds when we replace $G$ by a more general commutative algebraic group, which is the extension of an abelian variety by some commutative linear algebraic group $H$. However, it is easy to see that if $\dim H\ge 2$ and $H$ is not a torus (i.e., $G$ is not a semiabelian variety), then there are examples when an irreducible subvariety $V\subseteq G$ meets a finitely generated subgroup of $G(\bbk)$ in a Zariski dense subset and moreover $V$ is not a coset of an algebraic subgroup of $G$ (as predicted by the variant of the Mordell--Lang conjecture).

\begin{example}
\label{example:G_a}
If $G=\G_a^2$, then the graph of any polynomial of degree larger than $1$ with integer coefficients would contain infinitely many integral points, thus contradicting the corresponding Mordell--Lang principle for $G$.
\end{example}

\begin{example}
If $G=\G_a\times \G_m$, then the diagonal subvariety $\Delta\subseteq G$ contains infinitely many points of the subgroup $\Gamma$ spanned by $(1,1)$ and $(0,2)$; more precisely, $(2^n,2^n)\in \Delta$ for each $n\in\N$. However, $\Delta$ is not a coset of a subgroup of $G$.   
\end{example}  

Vojta~\cite[page~134]{Voj96} noted 
the example, similar to our 
Example~\ref{example:G_a}, of the subvariety
of $\G_a^2$ defined by Pell's equation as evidence for 
his ``doubt that this [Mordell--Lang] result can be 
generalized to a larger class of group varieties''.
We see, instead, that these examples suggest 
that the only commutative algebraic group $G$ for which a variant of the Mordell--Lang conjecture might hold is an extension of an abelian variety by a single copy of $\G_a$. Assuming the Bombieri--Lang conjecture\footnote{Bombieri formulated this conjecture only in the case of surfaces.}, we prove in Section~\ref{sec:number field} that a variant of the Mordell--Lang conjecture holds for an algebraic group isomorphic to a product of an abelian variety with the one-dimensional additive group. First we state the aforementioned conjecture of Bombieri--Lang; for more details on this famous conjecture, see \cite[Chapter~14]{BG06}. We also note that the Bombieri--Lang conjecture is a special case of Vojta's conjectures which play a central role in arithmetic geometry (see \cite[Conjecture~14.3.2~and~its~Remark~14.3.7]{BG06}).

\begin{conjecture}[Bombieri--Lang--Vojta]
\label{conj:BLV}
Let $X$ be a projective variety of general type defined over $\Qbar$. Then for each number field $\K$, the set $X(\K)$ is not Zariski dense in $X$. 
\end{conjecture}
A variety $X$ is of \emph{general type} if its canonical divisor $K_X$ is ample; often times Conjecture~\ref{conj:BLV} is formulated under the weaker assumption that $K_X$ is a big divisor. 
If $\dim X=1$, then Conjecture~\ref{conj:BLV} is equivalent to the well-known Mordell conjecture, proven by Faltings \cite{Fal83}. We prove the following result.

\begin{theorem}
\label{main theorem}
Let $A$ be an abelian variety defined over $\Qbar$ and let $\Gamma$ be a finitely generated subgroup of $(\G_a\times A)(\Qbar)$. If Conjecture~\ref{conj:BLV} holds, then for each  subvariety $V\subseteq \G_a\times A$, the intersection $V(\Qbar)\cap\Gamma$ is a finite union of cosets of subgroups of $\Gamma$.
\end{theorem}

The problem of proving unconditionally the results of this paper in the
most general case seems deeply linked with the conjecture of Vojta for a
generically finite cover of an abelian variety. In turn, this seems to
be out of reach of the present methods, unless one inserts additional
assumptions.

Unconditionally, we can prove Theorem~\ref{main theorem} when $\dim A = 1$ (see Theorem~\ref{thm:elliptic}), and also when $V$ is birational to a subvariety of an abelian variety (see Remark~\ref{rem:unconditional}). In the case when $V$ is birational to a subvariety of some abelian variety, the result follows as a consequence of Faltings' result \cite[Theorem~2]{Fal91} regarding the finiteness of the number of $S$-integral points on an abelian variety with respect to an ample divisor (see also Remark~\ref{rem:unconditional} and our proof of Theorem~\ref{main theorem}). In the case $A$ is an elliptic curve, we can prove a more general result (inspired by the results of \cite{CMZ}), valid for any commutative algebraic group of dimension $2$ that is an  extension of an elliptic curve by a copy of the additive group.

\begin{theorem}
\label{thm:elliptic}
Let $\K$ be a finitely generated field of characteristic $0$, let $E/\K$ be an elliptic curve and let $G$ be a commutative algebraic group which is an extension of $E$ by $\bG_a$, i.e., there is a short exact sequence of connected algebraic groups:
$$0\lra \bG_a\lra G\lra E\lra 0.$$
Let $\Gamma$ be a subgroup of $G(\K)$ such that $\Gamma\cap \bG_a$ is finitely generated, and let
$T$ be a subset of $\Gamma$.
Then the Zariski closure of $T$ is  a finite union of
translates of algebraic subgroups of $G$.
\end{theorem}

Assuming that $A$ is an abelian variety of $\bbk/\Qbar$-trace $0$ (i.e., there is no nonconstant morphism between $A$ and some abelian variety defined over $\Qbar$), then we can prove unconditionally the conclusion from Theorem~\ref{main theorem} even in the more general case when we intersect a subvariety $V\subseteq \G_a\times A$ with a finite rank subgroup.

\begin{theorem}
\label{thm:function field}
Let $\bbk$ be an algebraically closed field of characteristic $0$, let $A$ be an abelian variety such that ${\rm Tr}_{\bbk/\Qbar}(A)=0$. Then each subvariety $V\subseteq \G_a\times A$ intersects a finite rank subgroup $\Gamma\leq (\G_a\times A)(\bbk)$ in a finite union of cosets of subgroups of $\Gamma$.
\end{theorem}

We prove Theorem~\ref{thm:function field} in Section~\ref{sec:function field}  using the ideas introduced by Hrushovski \cite{Hrushovski} for his proof of the function field version of the Mordell--Lang conjecture. We note that Theorem~\ref{thm:function field} fails if the abelian variety were defined over a number field. Indeed, if $A$ is an abelian variety whose $\Q$-rational points are Zariski dense, while  $V\subseteq \G_a\times A$ is the graph of any non-constant rational function $f\colon A\lra \bP^1$ defined over $\Q$, and $\Gamma \coloneqq (\G_a\times A)(\Q)$, then $V(\Q)\cap\Gamma$ is Zariski dense in $V$, even though $V$ is not a coset of an algebraic subgroup of $\G_a\times A$.


\section{Proofs of Theorems~\ref{main theorem} and \ref{thm:elliptic}}
\label{sec:number field}

\begin{proof}[Proof of Theorem~\ref{main theorem}.] 
We only need to prove the case when $V$ is irreducible. Hence we may assume that $V$ is integral. Next, it suffices to prove that if $V\cap \Gamma$ is Zariski dense in $V$ then $V$ is a translate of an algebraic subgroup of $\G_a\times A$. Let $\pi_2\colon \G_a\times A\lra A$ be the usual projection morphism. Then by Faltings' theorem \cite{Fal94}, we may replace $A$ by the Zariski closure of $\pi_2(V\cap \Gamma)$ which is a translate of an abelian subvariety of $A$, and assume that $\pi_2|_V\colon V\lra A$ is dominant.
If $\pi_2|_V$ is of relative dimension $1$, then $V = \G_a \times A$ and we are done. 
So we only need to consider the case that $\pi_2|_V$ is generically finite. Let $\pi_1\colon \G_a\times A\lra \G_a$ be the usual projection morphism. If $\pi_1|_V\colon V\lra \G_a$ is not dominant, i.e., $\pi_1(V)$ is a point, then $V$ is contained in a fibre $F\cong A$ of $\pi_1$.
It follows that $V = F$ by the dimension reasoning and hence $V$ is a translate of the abelian variety $A$.
Hence, from now on, we may assume further that $\pi_1|_V\colon V\lra \G_a$ is dominant. We will prove that in this case, our hypotheses yield a contradiction.

Let $\bP^1\times A$ be the compactification of $\G_a\times A$ (such that $\bP^1 - \G_a = \{\infty\}$) and $\overline{V}$ the Zariski closure of $V$ in $\bP^1\times A$. Also, by a slight abuse of notation, we still denote by $\pi_2$ the projection morphism $\bP^1\times A\lra A$. 
Then by the Stein factorization, there exists a normal projective variety $Y$ endowed with a birational morphism $\iota\colon \overline{V}\lra Y$ and with a finite morphism $f\colon Y\lra A$ such that $\pi_2|_{\overline{V}}=f\circ \iota$. 

Let $\K$ be a number field such that $V$, $\overline{V}$, $Y$, $\iota$ and $f$ are all defined over $\K$, and also $\Gamma\leq (\G_a\times A)(\K)$. Since $V(\K)\cap \Gamma$ is Zariski dense in $V$, we conclude that $Y(\K)$ is Zariski dense in $Y$.

Applying Kawamata's structure theorem \cite[Theorem~13]{Kawamata} to the finite morphism $f\colon Y \lra A$, there exists a finite \'etale cover $\phi\colon \widetilde{Y}\lra Y$ such that $\widetilde{Y}$ is isomorphic to the direct product $\widetilde{B}\times W$, where $\widetilde{B}$ is a finite \'etale cover of an abelian subvariety $B$ of $A$ and $W$ is a normal projective variety of general type, i.e., $\kappa(W) = \dim W$. Since $Y(\K)$ is Zariski dense in $Y$ and $\phi\colon \widetilde{Y}\lra Y$ is \'etale, the Chevalley--Weil theorem (see \cite[p.~585]{C-Z} and \cite[Theorem~10.3.11]{BG06}) yields that there exists a finite extension $\L/\K$ such that $\widetilde{Y}(\L)$ is Zariski dense in $\widetilde{Y}$. In particular, at the expense of replacing $\L$ by another finite extension, we obtain that $W(\L)$ is Zariski dense in $W$, which contradicts Conjecture~\ref{conj:BLV}, if $\dim W > 0$. 

So, from now on, we may assume that $W$ is a point, which yields that $\widetilde{Y}$ and therefore $Y$ itself is an abelian variety. Since $Y$ is birational to $V$, composing this birational map with $\pi_1|_V$, we obtain a non-constant rational function $g\colon Y\lra \bP^1$ (note that $\pi_1|_V\colon V\lra \G_a$ is dominant). Since $\Gamma$ is a finitely generated subgroup of $(\G_a \times A)(\K)$, we have that $\pi_1(\Gamma)$ is a set of $S$-integral points in $\K$ with respect to a suitable finite set $S$ of places of $\K$.

We let $D \coloneqq g^*(\{\infty\})$ be the divisor of $Y$ which is the pullback of the point at infinity for the inclusion $\bG_a\subseteq \bP^1$; note that $D$ is a divisor because $g$ is non-constant. We let $\tilde{\iota}\colon V\lra Y$ be the corresponding birational map and then for each point $x\in \tilde{\iota}(\Gamma\cap V)$ (note that $\Gamma\cap V$ is Zariski dense in $V$), we have that $g(x)$ is $S$-integral with respect to the divisor $D$ of $Y$. If $D$ is ample, then Faltings' theorem \cite[Theorem~2]{Fal91} yields a contradiction to the fact that there exist infinitely many such $S$-integral points (for a variant of Faltings' theorem in the context of semiabelian varieties, see \cite{Voj99}). We show next that the general case reduces to this special case.

Assume now that $D$ is not ample. Let $C$ be the connected component of the stabilizer of $D$ in $Y$. Then $C$ is an abelian subvariety of $Y$ of positive dimension. Let $Z$ be a complement of $C$ in $Y$, i.e., $Z$ is a proper abelian subvariety of $Y$ such that $Y$ is isogenous to $C\times Z$. Therefore, without loss of generality, we may replace $Y$ by $C\times Z$. We let $h \coloneqq g|_Z$, which is still a non-constant rational function $Z\lra \bP^1$ and moreover, $h^*(\{\infty\})$ is an ample divisor of $Z$. Another application of Faltings' \cite[Theorem~2]{Fal91} provides a contradiction, which finishes our proof of Theorem~\ref{main theorem}.
\end{proof}

\begin{remark}
\label{rem:simple}
Using the notation as in Theorem~\ref{main theorem}, if $A$ is a {\it simple} abelian variety, then one does not need to use \cite[Theorem~13]{Kawamata} to finish the proof. Indeed, \cite{Kawamata} was employed only in the case the finite morphism $f\colon Y\lra A$ is ramified, because in the case $f$ is unramified, we immediately derive that $Y$ must itself be an abelian variety and proceed as in the proof of Theorem~\ref{main theorem} invoking only Faltings' theorem \cite{Fal91} regarding the finiteness of the number of $S$-integral points on an abelian variety with respect to an ample divisor. Now, if $f$ is ramified, the canonical divisor $K_Y$ of $Y$ is  the ramification divisor of $f$, and moreover since $f$ is a finite map, we obtain that $K_Y=f^*(D_f)$ for some effective divisor $D_f$ of $A$. If $A$ is simple, then each nontrivial effective divisor of it is ample and therefore we obtain that $K_Y$ is big, which still allows us to apply the Bombieri--Lang--Vojta conjecture to obtain a contradiction.
\end{remark}

\begin{remark}
\label{rem:unconditional}
As shown in our proof of Theorem~\ref{main theorem}, the only point in which we employed the validity of Conjecture~\ref{conj:BLV} is for the case when the finite morphism $Y\lra A$ is ramified. In particular, this means that Theorem~\ref{main theorem} holds unconditionally if the subvariety $V\subseteq \G_a\times A$ is birational to an abelian variety. Furthermore, if $V\subseteq \G_a\times A$ is birational to a subvariety $Y$ of some arbitrary abelian variety $B$, then the assumption that $V$ contains a Zariski dense set of $\K$-rational points yields that $Y$ contains a Zariski dense set of rational points and therefore, Faltings' theorem \cite{Fal94} yields that $Y$ must be a coset of an abelian subvariety of $B$. So, $V$ is birational to an abelian variety itself and we are done using Faltings' theorem \cite{Fal91} regarding the finiteness of the $S$-integral points on an abelian variety.
\end{remark}

For the general case of a non-split extension $G$ of an arbitrary abelian variety $A$ (defined over a field of characteristic $0$) by a copy of the additive group, the corresponding variant of the Mordell--Lang conjecture is quite subtle. However, we can settle unconditionally  the case when $A$ is an elliptic
curve.

\begin{proof}[Proof of Theorem~\ref{thm:elliptic}.]
First, we observe that since $\Gamma$ projects to $E(\K)$, which is a finitely generated group (due to the classical Mordell--Weil theorem), we get that $\Gamma$ must itself be finitely generated (since so is its intersection with $\bG_a(\K)$ by the assumption). Hence, our goal is to show that if $V\subseteq G$ is an irreducible curve with the property that $V\cap \Gamma$ is Zariski dense in $V$, then $V$ must be a coset of a one-dimensional algebraic subgroup of $G$. We have two cases: either $G$ is a split extension, or not.

{\it Case 1.} $G$ is a split extension. So, at the expense of replacing $G$ by an isogenous copy of it and also replace $\bbk$ by a finite extension, we may assume $G=\bG_a\times E$.

If $V$ does not project dominantly onto one of the two factors of $\bG_a\times E$, then we obtain the desired conclusion. So, assume the curve $V\subseteq \G_a\times E$ projects dominantly onto both factors of $\G_a\times E$; using the hypothesis that $V$ contains infinitely many points of the subgroup $\Gamma\leq (\G_a\times E)(\K)$, then we derive a contradiction. We observe that $V$ must have positive genus since it projects dominantly onto the elliptic curve $E$; then we derive (similar to the proof of Theorem~\ref{main theorem}) a contradiction due to the finiteness of the number of $S$-integral points on a curve of positive genus.
 
{\it Case 2.} $G$ is a non-split extension.

We settle this case using the fact that $G$ does not contain complete
curves (see \cite{CMZ} for this and other facts on such group extensions). So,
if the curve $V\subseteq G$ contained infinitely many
points from a finitely generated subgroup $\Gamma$ of $G$, then these points
would be $S$-integral with respect to the complement of $V$ in the projective
closure of it, where $S$ is a finite set of places containing those of bad
reduction either for $G$ or for a set of generators for $\Gamma$. (See also \cite{CMZ} for an explicit projective embedding of $G$, obtained
first by Serre, in a letter to Masser reproduced in the aforementioned paper.)   
But since $V$ is affine, Siegel's Theorem would entail that $V$ has genus $0$;
hence $V$ could not dominate $E$ and would be a translate of $\G_a$, as required. 
\end{proof}


\section{Proof of Theorem~\ref{thm:function field}}
\label{sec:function field}

\noindent
We work with the differential algebraic methods of~\cite{Bui92,Hrushovski}
and recommend the book~\cite{MMP} for general background.

We begin by endowing $\bbk$ with a derivation $\partial \colon \bbk \to \bbk$ for which 
the field of constants $\bbk^{\partial} \coloneqq \{ x \in \bbk ~:~ \partial(x) 
= 0 \}$ is the field $\Qbar$ of algebraic numbers.  Since the 
statement of Theorem~\ref{thm:function field} becomes only formally 
more difficult with $\bbk$ replaced by a larger field, we may, and do, 
replace $\bbk$ with its differential closure, which is still an 
algebraically closed differential field having field of constants 
equal to $\Qbar$.  For the sake of readability, we shall 
identify algebraic (and differential algebraic) varieties 
with their sets of $\bbk$-points.

By the usual reductions, we may assume that $V$ is irreducible, 
contains the identity element of the group, and has a trivial 
stabilizer.  We are charged with showing that if $V\cap \Gamma$ is Zariski dense in $V$ (where $\Gamma\leq (\bG_a\times A)(\bbk)$ is a group of finite rank), then $V$ must consist of a
single point.   

We find a differential algebraic subgroup $\widetilde{\Gamma} \leq 
\G_a\times A$ for which $\Gamma \leq \widetilde{\Gamma}$ and 
$\widetilde{\Gamma}$ has finite Morley rank (see~\cite[(5.1)]{Mar00}). 
Let $\Xi$ be the image of $\widetilde{\Gamma}$ under the projection 
map $ \G_a\times A \to A$ and let $B$ be the connected component of the 
identity of the Zariski closure of $\Xi$.  Since it is an algebraic 
subgroup of $A$, then $B$ also has $\bbk/\Qbar$-trace zero.   Hence, 
by~\cite[Proposition 2.6]{HrSo} (also proven in~\cite{PiZi03})
the Manin kernel $B^\sharp$ of $B$ is locally modular and is thus
orthogonal to the field of constants.  By~\cite[Lemma~4.2]{Pil96}, every
Zariski dense differential algebraic subgroup of $B$ contains $B^\sharp$; 
thus, $B^\sharp \leq \Xi$.  Consider also the image $\Upsilon$ of 
$\widetilde{\Gamma}$ under the projection $\G_a\times A \to \G_a$.  
Since $\Upsilon$ is a finite Morley rank subgroup of the additive 
group, it is a finite dimensional vector space over the field of 
constants and is, therefore, fully orthogonal to $B^\sharp$.  The group 
$\widetilde{\Gamma} \cap (\Upsilon\times B^\sharp)$ is a differential
algebraic subgroup of $\Upsilon\times B^\sharp$ which projects 
onto $B^\sharp$.  By the orthogonality of $B^\sharp$ and $\Upsilon$, 
every differential algebraic subvariety is a union of products $T \times S$
where $S \subseteq B^\sharp$ and $T \subseteq \Upsilon$.  It follows that 
$\{0\}\times B^\sharp  \leq \widetilde{\Gamma}$.   

Let $Y \coloneqq \widetilde{\Gamma} \cap V$.   Since $\Gamma \cap V$ is 
Zariski dense in $V$, we have that the differential algebraic 
variety $Y$ is Zariski dense in $V$.  Since $V$ is irreducible as 
an algebraic variety, there exists some component $X$ of $Y$ which is 
Zariski dense in $V$.   Translating, we may assume that 
$X$ contains the identity element. 

By~\cite[Proposition 4.4]{Hrushovski} there is a differential algebraic
groups $H \leq \widetilde{\Gamma}$ for which $X$ is a union of cosets of
the connected component of $(H \cap (\{0\}\times B^\sharp))$ and 
$X$ is contained in $H + (\{0\}\times B^\sharp)$.  The group 
$(H \cap (\{0\}\times B^\sharp))$ is contained in the stabilizer of
$X$ and the Zariski closure of this group is contained in the 
stabilizer of $V$.   As we have reduced to the case that 
$V$ has a trivial stabilizer, $H \cap (\{0\}\times B^\sharp )$ is 
itself trivial.   

Using again that $B^\sharp$ is orthogonal to the 
field of constants, it follows that $H \leq \G_a\times \{ 0 \}$.  
Indeed, as before we let $\Phi$ be the projection of $H$ to $A$ and 
$\Psi$ be the projection of $H$ to $\G_a$.  As $\Psi \perp \Phi$
and $H \leq \Psi\times \Phi$, it must be that $H = \Psi \times \Phi$.
By~\cite[Lemma~4.2]{Pil96} again, $\Phi$ contains the 
Manin kernel $C^\sharp$ of the connected component $C$ 
of its Zariski closure, and $C^\sharp \leq B^\sharp$.  We know 
the group $(\{0\}\times B^\sharp ) \cap H$ contains $\{0\}\times C^\sharp $
and is trivial.  Hence $C^\sharp=\{0\}$ and therefore $C = \{ 0 \}$  (note that $C^\sharp$ is Zariski dense in $C$);  so, $H \leq \G_a\times \{ 0 \}$.  

Therefore, we know that $X \subseteq H + (\{0\}\times B^\sharp)$ and that 
the two groups in the sum are orthogonal.  Hence, 
$X$ may be expressed as $S + T$ with $S \subseteq H$ and $T 
\subseteq (\{0\}\times B^\sharp)$.  Since $B^\sharp$ is 
locally modular, by~\cite{HrPi} the set $T$ is a translate of a 
subgroup.  Using that $X$ has a trivial stabilizer, it follows that 
$T$ is a single point.  Taking Zariski closures, we see that 
$V$ is a translate of a subvariety of $\G_a\times \{ 0 \}$.  
Because $V$ has trivial stabilizer, this subvariety cannot be all of 
$\G_a\times \{ 0 \}$.  Thus, $V$ is a single point as we needed to show. This concludes our proof of Theorem~\ref{thm:function field}.

\begin{remark}
We expect that these differential algebraic 
techniques could be pushed to prove a relative Mordell--Lang
theorem for general commutative algebraic groups.   The statement we expect to 
be 	true is the following.  Let $G$ be a 
commutative algebraic group over the algebraically
closed field $\bbk$ of characteristic zero.  
Let $\Gamma \leq G(\bbk)$ be a finite rank 
subgroup and let $V \subseteq G$ be an irreducible
subvariety for which $\Gamma \cap V(\bbk)$ is 
Zariski dense in $V$.   Then there should be 
an algebraic subgroup $H \leq G$ of $G$, an 
algebraic group $J$ defined over $\Qbar$, 
an algebraic variety $Y \subseteq J$ also defined
over $\Qbar$, a point $a \in G(\bbk)$, and 
a map of algebraic groups $h\colon H \to J_\bbk$ (where
$J_\bbk$ is the base change of $J$ from $\Qbar$ 
to $\bbk$) so that $V = a + h^{-1} (Y_\bbk)$.   
When $G$ is an extension of an abelian variety 
of $\Qbar$-trace zero by $\G_a$, then this
conjecture would imply that $V$ must be a 
translate of an algebraic subgroup of $G$.  
What is missing from the existing literature on the
structure of differential algebraic groups is an 
appropriate analogue of the Socle Theorem~(\cite[Propositions 4.3 and 4.4]{Hrushovski}) for groups with non-rigid semipluriminimal 
socles. 
\end{remark}
 

\end{document}